\documentclass[12pt]{article}
\usepackage{amsthm,amsmath,amssymb}
\def\bbe{\mathbb E}

\def\bbn{\mathbb N}
\def\udel{\Delta}
\def\ugam{\Gamma}
\def\lam{\lambda}
\def\calD{{\cal D}}
\theoremstyle{plain}
\newtheorem{thm}{Theorem}
\theoremstyle{definition}
\newtheorem{rem}{Remark}

\begin{document}

\title{On a Generalization of the Arcsine Law}

\author{Christian Houdr\'e\thanks{Georgia Institute of Technology, 
School of Mathematics, Atlanta, GA 30332-0160; houdre@math.gatech.edu. 
This work was supported in part by the grants \#246283 and \#524678  from 
the Simons Foundation and by a Simons Foundation Fellowship, grant \#267336, while this author was visiting the 
LPMA of the Universit\'e Pierre et Maris Curie.} 
\and Trevis J. 
Litherland\thanks{LexisNexis, Inc., Alpharetta, GA 30005; trevis.litherland@lexisnexisrisk.com}}

\maketitle

\vspace{0.5cm}

\begin{abstract}
The gaps between the times, in a Weyl chamber, at which the sum of the increments of 
independent Brownian motions attains its maximum has a Dirichlet distribution.  

\end{abstract}



\vspace{0.5cm}


The Brownian functional $D_m$, $m \in \bbn$, introduced by Glynn and Whitt 
\cite{GW}, in the context of a queuing problem, is defined as:  
\begin{equation}\label{eq1.1}
D_m=\max_{{\scriptstyle 0=t_0\le t_1\le\cdots\le t_{m-1}\le t_m=1}}
\sum^m_{i=1} [B^{i}(t_i)-B^{i} (t_{i-1})],
\end{equation}
where $(B^{1}(t),\dots, B^{m}(t))$ is an $m$-dimensional standard Brownian
motion.  The asymptotic behavior, for $m\to +\infty$,  of this functional is well understood, in particular, 
since Baryshnikov \cite{Ba} as well 
as Gravner, Tracy and Widom \cite{GTW} identified its law as that of the law of the maximal eigenvalue of 
an $m\times m$ element of the Gaussian Unitary Ensemble (GUE).    (In fact, 
\cite{Ba} showed that the whole sequence  
$(D_m)_{m\ge 1}$ is identical in law to the sequence $(\lam_{m}^{*})_{m\ge 1}$,
where each $\lam^*_m$ is the largest eigenvalue of the $m \times m$ 
principal minor of an infinite GUE matrix.   Moreover, there is equality in law between 
all the spectra of all the principal  minors of an element of the GUE and some 
multidimensional Brownian functionals, 
naturally appearing as limits in subsequence problems; see \cite{BGH} and the references therein.)

The purpose of this note is to describe the distribution of the  (a.s.~unique) 
parameters $(\theta_1,\theta_2 ,\dots, \theta_{m-1})$ 
which maximize the functional $D_m$. Our result is the following generalization of the arcsine
law:
                                     ~
\begin{thm}\label{thm1.1}
 The density $f_m$, $m \ge 2$, of the parameters  $\tilde\theta^{(m)}
 := (\theta_1,\theta_2,\dots$, $\theta_{m-1})$,
$0=\theta_0\le\theta_1 \le \theta_2\le\cdots\le \theta_{m-1} \le \theta_m=1$, which maximize the
functional $D_m$ is given by
\begin{equation}\label{eq1.2}
f_m(\tilde\theta^{(m)})=\frac{\ugam\left(\frac m2\right)}{\pi^{m/2}}\, 
\theta_1^{-1/2}(1-\theta_{m-1})^{-1/2}
\prod^{m-1}_{i=2} (\theta_i-\theta_{i-1})^{-1/2}.
\end{equation}
That is, the $m$-vector of gap lengths 
$(\udel_1,\udel_2 ,\dots,\udel_m)$, where $\udel_i = \theta_i - \theta_{i-1}$,
$1\le i\le m$, 
has a Dirichlet distribution with parameters 
$(1/2, 1/2, \dots , 1/2)$.
\end{thm}

\begin{rem}\label{rem1.1}
 The gap distribution $\calD(1/2, 1/2,\dots , 1/2)$ shows that the 
optimal parameters are essentially evenly spaced, in the basic sense that 
$\bbe \theta_i = i/m$, for 
$1 \le i\le  m - 1$, as well as in the more universal sense that the
empirical distribution function of  $\tilde\theta^{(m)}$
 tends strongly towards the uniform
distribution on $[0, 1]$ as $m \to\infty$. 

\end{rem}
\begin{rem}\label{rem1.2}
Consecutive parameter values do not exhibit a strong repulsion-type 
property, unlike the behavior observed, for example, in GOE/GUE/GSE eigenvalues.  
The factors $(\theta_i - \theta_{i-1})^{-1/2}$ in \eqref{eq1.2}, 
in fact, introduce the bimodal phenomenon of a 
Beta$(1/2,  1/2)$ distribution, whose quartiles are found at 
$(2 -\sqrt 2)/4\approx 0.146$, $0.5$,
and $(2 +\sqrt2)/4\approx 0.854$. 
In effect, the gaps thus tend to be either quite small or
 large with respect to the natural length scale $1/m$.
\end{rem}

\begin{rem}	\label{rem1.3}
In his paper on the joint distribution of the location of the maximum of a Brownian motion 
and related quantities, Shepp \cite{Sh} provides a simple argument proving the a.s.~uniqueness of the
location parameter, an argument that generalizes to the a.s.~uniqueness of
$(\theta_1,\theta_2,\dots$, $\theta_{m-1})$.
\end{rem}

\begin{proof}
 We first note, and show, that that the parameters enjoy a time-reversal symmetry described by
\begin{equation}\label{eq1.3}
(\theta_1,\theta_2,\dots, \theta_{m-1})\buildrel d\over =(1-\theta_{m-1},
1-\theta_{m-2},\dots, 1-\theta_1).
\end{equation}

Indeed, to prove \eqref{eq1.3}, observe that for any fixed 
$0 = t_0\le t_1 \le\cdots\le t_{m-1}\le t_m = 1$, we have
\begin{align}\label{eq1.4}
\sum^m_{i=1} [B^i(t_i)-B^i(t_{i-1})]&\buildrel d\over =
\sum^m_{i=1} [B^i(1-t_{i-1})-B^i(1-t_i)]\nonumber\\
&\buildrel d\over = \sum^m_{i=1} [B^{m-i+1}(1-t_{i-1})-B^{m-i+1}(1-t_i)]\nonumber\\
&\buildrel d\over = \sum^m_{i=1} [B^i(1-t_{m-i})-B^i(1-t_{m-i+1})],
\end{align}
where the first two distributional equalities follow respectively
from the stationarity of the increments and from the 
independence of the coordinates of standard 
multidimensional Brownian motion. Thus, \eqref{eq1.3}
clearly follows.

Next, we note that since
\begin{align*}
D_2&=\max_{0\le t_1\le 1} [B^1(t_1)-B^1(0) + B^2(1)-B^2(t_1)]\\
&=B^2(1)+\max_{0\le t_1\le 1} [B^1(t_1)-B^2(t_1)],
\end{align*}
the problem reduces to finding the location $\theta_1$ of the maximum of a single
standard Brownian motion over $[0, 1]$. The solution to this problem, 
given by one of Paul L\'evy's arcsine law results, is well-known 
(see, e.g., Shepp's more general result in \cite{Sh}), 
and has density $f_2$ given by
\begin{equation}\label{eq1.5}
f_2(\theta_1)=\frac1\pi \frac1{\sqrt{\theta_1(1-\theta_1)}}\,,
\quad 0\le\theta_1\le 1,
\end{equation}
i.e., $\theta_1$ follows a Beta distribution with parameters 
$1/2$ and $1/2$.  Equivalently, 
$(\udel_1 ,\udel_2) := (\theta_1 , 1 - \theta_1 )$
 exhibits a $\calD(1/2, 1/2)$ distribution.

To complete the proof, we proceed by induction on $m$. 
Assume that the result is true for
$2, \dots, m$. Then, conditioning on $\theta_m$, 
and using scaling and the induction
hypothesis, we see that
\begin{align}\label{eq1.6}
f_{m+1}(\theta_1,\dots, \theta_m)&=f_{m+1}(\theta_1,\dots, \theta_{m-1}\mid
\theta_m)g_{m+1}(\theta_m)\nonumber\\
&=\left\{f_m\left(\frac{\theta_1}{\theta_m} ,\dots, 
\frac{\theta_{m-1}}{\theta_m}
\right)\frac1{(\theta_m)^{m-1}}\right\}g_{m+1}(\theta_m),
\end{align}
where $g_{m+1}$ is the marginal density of $\theta_m$. 
Similarly, conditioning this time
on $\theta_1$, and again using scaling and the induction hypothesis, 
we find that
\begin{align}\label{eq1.7}
f_{m+1}&(\theta_1,\dots, \theta_m)=f_{m+1}(\theta_2,\dots,\theta_m\mid
\theta_1)g_{m+1}(1-\theta_1)\nonumber\\
&=\left\{f_m\left(\frac{\theta_2-\theta_1}{1-\theta_1} ,\dots, 
\frac{\theta_m-\theta_1}{1-\theta_1}
\right)\frac1{(1-\theta_1)^{m-1}}\right\}g_{m+1}(1-\theta_1),
\end{align}
where we have also made use of the time-reversal symmetry property  
\eqref{eq1.3}. Equating
\eqref{eq1.6} and \eqref{eq1.7}, we can solve for $g_{m+1}$ as follows:
\begin{align}\label{eq1.8}
&f_m\left(\frac{\theta_1}{\theta_m} ,\dots, 
\frac{\theta_{m-1}}{\theta_m} \right)
\frac1{(\theta_m)^{m-1}}\, g_{m+1}(\theta_m)\nonumber\\
&\quad =\frac{\ugam \left(\frac m2\right)}{\pi^{m/2}}
\left(\frac{\theta_1}{\theta_m}\right)^{-1/2}
\left(\frac{\theta_m-\theta_{m-1}}{\theta_m}\right)^{-1/2}
\prod^{m-1}_{i=2}\left(\frac{\theta_i-\theta_{i-1}}{\theta_m}\right)^{-1/2}
\frac{g_{m+1}(\theta_m)}{(\theta_m)^{m-1}}\nonumber\\
&\quad=\frac{\ugam\left(\frac m2\right)}{\pi^{m/2}}\,
(\theta_1)^{-1/2}
\prod^m_{i=2} (\theta_i-\theta_{i-1})^{-1/2}
\frac{g_{m+1}(\theta_m)}{(\theta_m)^{\frac m2-1}}\nonumber\\
&\quad = f_m\left(\frac{\theta_2-\theta_1}{1-\theta_1},\dots, 
\frac{\theta_m-\theta_1}{1-\theta_1}\right)\frac1{(1-\theta_1)^{m-1}}\,
g_{m+1}(1-\theta_1)\nonumber\\
&\quad = \frac{\ugam \left(\frac m2\right)}{\pi^{m/2}}
\left(\frac{\theta_2-\theta_1}{1-\theta_1}\right)^{-1/2}
\left(\frac{1-\theta_m}{1-\theta_1}\right)^{-1/2}
\prod^m_{i=3}\left(\frac{\theta_i-\theta_{i-1}}{1-\theta_1}\right)^{-1/2}
\frac{g_{m+1}(1-\theta_1)}
{(1-\theta_1)^{m-1}}\nonumber\\
&\quad =\frac{\ugam\left(\frac m2\right)}{\pi^{m/2}}
(1-\theta_m)^{-1/2}
\prod^m_{i=2}(\theta_i-\theta_{i-1})^{-1/2}
\frac{g_{m+1}(1-\theta_1)}{(1-\theta_1)^{\frac m2-1}}\,,
\end{align}
and, after cancellation of terms, we have
$$\frac{g_{m+1}(\theta_m)}{(1-\theta_m)^{\frac12-1}(\theta_m)^{\frac m2-1}}=
\frac{g_{m+1}(1-\theta_1)}{\theta^{\frac12-1}_1(1-\theta_1)^{\frac m2-1}}\,,
$$
for all $0\le \theta_1\le\theta_m \le1$. 
But this implies that for all $0\le t\le 1$,
$$g_{m+1}(t)\propto t^{\frac m2-1}(1-t)^{-\frac12},$$
and hence, after normalization, we have
\begin{equation}\label{eq1.9}
g_{m+1}(t)=\frac{\ugam \left(\frac{m+1}2\right)}{\ugam\left(\frac m2\right)
\ugam \left(\frac12\right)}\, t^{\frac m2-1}(1-t)^{-\frac12},
\end{equation}
the density of a Beta$(m/2, 1/2)$ distribution. Plugging 
\eqref{eq1.9} back into \eqref{eq1.6}, and recalling that
$\Gamma(1/2) = \sqrt{\pi}$, we find that
\begin{align}\label{eq1.10}
f_{m+1}(\theta_1,\dots,\theta_m)&= f_m\left(\frac{\theta_1}{\theta_m},\dots,
\frac{\theta_{m-1}}{\theta_m}\right)\frac1{(\theta_m)^{m-1}}\,
g_{m+1}(\theta_m)\nonumber\\
&=\frac{\ugam \left(\frac m2\right)}{\pi^{m/2}}\,
\theta_1^{-1/2}
\prod^m_{i=2}(\theta_i-\theta_{i-1})^{-1/2}
\frac1{(\theta_m)^{\frac m2-1}}\nonumber\\
&\quad\times \frac{\ugam \left(\frac{m+1}2\right)}{\ugam \left(\frac m2\right)
\ugam \left(\frac12\right)}\,(\theta_m)^{\frac m2-1} (1-\theta_m)^{-1/2}
\nonumber\\
&= \frac{\ugam\left(\frac{m+1}2\right)}{\pi^{(m+1)/2}}\,
\theta_1^{-1/2}(1-\theta_m)^{-1/2}
\prod^m_{i=2} (\theta_i-\theta_{i-1})^{-1/2},
\end{align}
which proves the result.
\end{proof}

This above result suggests further explorations of other properties of closely related Brownian functionals.  

Shepp \cite{Sh} derived the {\it joint} density of 
the location $\theta$ of the maximum, the maximum $M_T$, and its final 
value $B(T)$ of a single Brownian motion over a finite
interval $[0, T ]$. Similarly, in our problem, it would be of
interest to determine the joint distribution of the location parameters  
$\tilde\theta^{(m)}$, $D_m$, and the terminal values $(B^1(1),...,B^m(1)).$

As already mentioned, \cite{BGH} obtains maximal Brownian functionals representations for the spectra 
of all the principal minors of a GUE matrix.    It would also be of interest to determine the law of 
the locations where these 
maxima are attained.

Another interesting question is the following. Given a random sample of
parameter locations $\{\tilde\theta^{(m,j)}\}^n_{j=1}$
from a $\calD(1/2, \dots,1/2)$ distribution, how do
$D_m$ and its empirical counterpart
$$D^n_m:=\max_{1\le j\le n}\sum^m_{i=1} \left[B^i
\left(\tilde\theta^{(m,j)}_{i}\right)-
B^i\left(\tilde\theta^{(m,j)}_{i-1}\right)\right],$$
compare (in expectation, asymptotically, etc.)? 

Finally, instead of defining
$D^n_m$ in terms of random samples 
$\{\tilde\theta^{(m,j)}\}^n_{j=1}$, one might ask what deterministic 
values of $\{\tilde\theta^{(m,j)}\}^n_{j=1}$ would maximize 
$\bbe D^n_m$, an idea analogous to the
numerical integration task of finding points for performing Gaussian 
quadrature.

\end{document}